\documentclass{amsart} 
\newtheorem{theorem}{Theorem}[section]
\newtheorem{proposition}[theorem]{Proposition}
\newtheorem{definition}[theorem]{Definition}
\newtheorem{lemma}[theorem]{Lemma}
\newtheorem{remark}[theorem]{Remark}
\newtheorem{corollary}[theorem]{Corollary}

\newcommand{\Br}{\mathrm{Br}}
\newcommand{\Z}{\mathbb{Z}}
\newcommand{\N}{\mathbb{N}}
\newcommand{\twist}{\mathfrak{I}(\theta)}

\newcommand{\id}{\mathrm{id}}
\newcommand{\fix}{\mathrm{Fix}(\theta)}

\begin{document}

\title{Fixed points of zircon automorphisms} 

\author{Axel Hultman}
\address{Department of Mathematics, KTH, SE-100 44, Stockholm, Sweden}
\email{axel@math.kth.se} 
\begin{abstract}
A {\em zircon} is a poset in which every principal
order ideal is finite and equipped with a so-called special
matching. We prove that the subposet 
induced by the fixed points of any automorphism of a zircon is itself
a zircon. This provides a natural context in which to view recent results on
Bruhat orders on twisted involutions in Coxeter groups.
\end{abstract} 
\maketitle

\section{Background and results}
Let $P$ be a partially ordered set ({\em poset}). A {\em matching} on
$P$ is an involution $M:P\to P$ such that $M(p)\lhd p$ or $p\lhd M(p)$
for all $p\in P$, where $\lhd$ denotes the covering relation of
$P$. In other words, $M$ is a graph-theoretic (complete) matching on
the Hasse diagram of $P$. 
\begin{definition}
Suppose $M$ is a matching on a poset $P$. Then, $M$ is called {\em
  special} if for all $p,q\in P$ with $p\lhd q$, we either have
  $M(p)=q$ or $M(p)<M(q)$. 
\end{definition}

The term ``special matching'' was coined by Brenti \cite{brenti}. In
the context of an Eulerian poset, a special matching is another way to
think of a {\em compression labelling} as defined by du Cloux \cite{ducloux}.

\begin{definition}\label{de:zircon}
A poset $P$ is a {\em zircon} if for any non-minimal element $x\in P$,
the subposet induced by the principal order ideal $\{p\in P\mid p\leq x\}$ is
finite and has a special matching. 
\end{definition}

Zircons were defined by Marietti in \cite{marietti}. Actually, his
definition differs somewhat from ours, but Proposition
\ref{pr:equal} below shows that they are equivalent. We have chosen
to use our definition because it is typically more convenient to check the
finiteness condition in Definition \ref{de:zircon} rather than finding
a rank function as required by the definition in \cite{marietti}.

The motivation to introduce zircons comes from the fact that they
mimic the behaviour of Coxeter groups ordered by the Bruhat
order. More precisely, the Bruhat order ideal below a non-identity
element $w$ in a Coxeter group has a special matching given by
multiplication with any descent of $w$. The finiteness condition in
Definition \ref{de:zircon} is trivially satisfied, implying that the
Bruhat order on any Coxeter group is a zircon.

In this note, we study the fixed points of automorphisms of
zircons. Our main results are the next theorem and its corollary. The
proofs are postponed to Section \ref{se:proofs}.

Say that a poset is {\em bounded} if it has unique maximal
and minimal elements.

\begin{theorem}\label{th:main}
Suppose $P$ is a finite, bounded poset equipped with a special
matching $M$. Let $\phi$ be an automorphism of $P$. Then, the subposet
of $P$ induced by the fixed points of $\phi$ has a special matching.
\end{theorem}

\begin{corollary}\label{co:main}
The fixed points of any automorphism of a zircon induce a subposet
which is itself a zircon.
\end{corollary}

Now, we briefly describe our reasons for being interested in results of this
kind. We refer to \cite{humphreys} or \cite{BB} for a thorough account
of the theory of Coxeter groups and their Bruhat orders. 

Let $(W,S)$ be a finitely generated Coxeter system. For $X\subseteq
W$, let $\Br(X)$ denote the subposet of the Bruhat order on $W$ which
is induced by $X$. A fundamental result due to Bj\"orner and Wachs \cite{BW}
asserts that the (order complex\footnote{The order complex of a poset
  is the (abstract) simplicial complex whose simplices are the 
totally ordered subsets.} of) the open intervals in
$\Br(W)$ are homeomorphic to spheres.  

An interesting subposet of $\Br(W)$ is induced by the involutions in
$W$. More generally, if $\theta:W\to W$ is an involutive group
automorphism which preserves the generating set $S$, the set of {\em twisted
involutions} is
\[
\twist = \{w\in W\mid \theta(w) = w^{-1}\}.
\]
The ordinary involutions are obtained by taking $\theta$ to be the
trivial automorphism. Richardson and Springer \cite{RS,RS2} showed
 that $\Br(\twist)$ is of importance to the study of cell
 decompositions of certain symmetric varieties. 

In \cite{hultman2} it was shown that $\Br(\twist)$, just as $\Br(W)$, has
the property that every open interval is a sphere. The method of proof was to
show that $\Br(\twist)$, too, is a zircon (although this terminology
was not used), and then observing that a result of Dyer \cite{dyer}
implies that every open interval in any zircon is homeomorphic to a sphere. 

Earlier, another approach to the topology of $\Br(\twist)$ was followed in
\cite{hultman} where it was observed that $\Br(\twist)$ is the
subposet of $\Br(W)$ induced by the fixed points of the involutive poset
automorphism given by $w \mapsto \theta(w^{-1})$. Invoking Smith
theory on automorphisms of spheres \cite{smith}, this permitted the
conclusion that the intervals in $\Br(\twist)$ are homology spheres
over the integers modulo $2$. 

To summarise, we have a zircon (namely $\Br(W)$) whose intervals are 
spheres. We construct an involution on it whose induced subposet of
fixed points (namely $\Br(\twist)$) also turns out to form a
zircon. Therefore, the intervals in this fixed point poset are not
only $\Z_2$-homology spheres (as implied by Smith theory) but actual
spheres.

Corollary \ref{co:main} explains this behaviour by showing
that, in fact, any automorphism of any zircon has a zircon as induced
fixed point poset.

\begin{remark}{\em 
Let $\theta$ and $W$ be as above. Then, $\theta$ is a poset
  automorphism of $\Br(W)$. Moreover, it is known \cite{hee, muhlherr,
  steinberg} that the fixed points $\fix$ themselves form a Coxeter
  group. It was shown in \cite{hultman}, and independently by Nanba
  \cite{nanba}, that the subposet of $\Br(W)$ induced by $\fix$
  coincides with the Bruhat order on $\fix$. In particular, this is
  another situation where the fixed points of a zircon automorphism
  are known to form a zircon.}
\end{remark}

\section{Proofs}\label{se:proofs}

The next lemma provides an extremely useful tool when dealing with posets
with special matchings. For Bruhat orders it was established by
Deodhar \cite[Theorem 1.1]{deodhar}. Brenti \cite[Lemma 4.2]{brenti2}
proved the general case under the assumption that $P$ is
graded. This assumption is, however, not essential to the proof. For
the reader's convenience, we restate Brenti's proof with slight adjustments.

A poset is {\em locally finite} if every interval is finite.

\begin{lemma}[Lifting property]
Suppose $P$ is a locally finite poset with a special matching
$M$. Choose $x,y\in P$ with $x < y$ and $M(y)<y$. Then,
\begin{enumerate}
\item[(i)] $M(x)\leq y$.
\item[(ii)] $M(x)<x \Rightarrow M(x) < M(y)$. 
\end{enumerate}
\begin{proof}
We proceed by induction over the length of a shortest non-refinable
chain $c$ in $[x,y]$, both statements following directly from the
definition of special matchings if $x\lhd y$.

For the first assertion, we may assume $M(x) > x$. Choose $z\in c$
such that $x\lhd z$. If $M(x) = z$, we are done. Otherwise, we have
$M(z) > z, M(x)$ since $M$ is a special matching. By the induction
assumption, $M(z) \leq y$, so that $M(x) \leq y$.

To prove the second claim, suppose $M(x)<x$. Pick $z\in c$ with $z\lhd
y$. In case $M(y) = z$, there is nothing to prove. Otherwise, because
$M$ is special, $M(z)<M(y)$. By induction, $M(x) < M(z)$, and the proof
is complete. 
\end{proof}
\end{lemma}

We are now ready to prove the main results.

\begin{proof}[Proof of Theorem \ref{th:main}]
The automorphism $\phi$ is of finite order $N$ since it is a permutation
of a finite set. Each automorphism $\phi^k$, $k\in [N] = \{1, \dots, N\}$,
transforms $M$ into a special matching $M_k$ on $P$. In particular, $M
= M_N$.

Given $p\in P$, let  
\[
C(p) = \{q \in P\mid q = M_{i_t}\circ M_{i_{t-1}}\circ \dots \circ
M_{i_1}(p)\text{ for some }i_1,\dots,i_t\in [N]\}.
\]
In other words, $C(p)$ consists of the elements in the same connected
component as $p$ in the graph we obtain from the Hasse diagram of $P$
by throwing away the edges that are not used by any of the matchings
$M_k$. By abuse of notation, we also let $C(p)$ denote the subposet of
$P$ induced by this set. 

Given $q\in C(p)$, we may write $q = M_{i_t} \circ \dots \circ
M_{i_1}(p)$ for suitably chosen $i_j\in [N]$. Now define $q^\prime\in
C(p)$ by $q^\prime = M^\prime_{i_t} \circ \dots \circ
M^\prime _{i_1}(p)$, where we recursively have defined
\[
M^\prime_{i_j} = 
\begin{cases}
M_{i_j} & \text{if } M_{i_j}\circ M^\prime_{i_{j-1}} \dots \circ
  M^\prime_{i_1}(p) < M^\prime_{i_{j-1}}\circ \dots \circ M^\prime_{i_1}(p),\\
\id & \text{otherwise.}
\end{cases}
\]
Repeatedly applying the lifting property, we obtain $M_{i_j}\circ
\dots \circ M_{i_1}(p)\geq M^\prime_{i_j}\circ
\dots \circ M^\prime_{i_1}(p)$ for all $j\in [t]$. In particular,
$q^\prime \leq q$. 

By construction, $p \geq q^\prime$. Furthermore, if $q$ is a minimal
element in $C(p)$, we have $q = q^\prime$. Thus, $C(p)$ has a unique
minimal element. A completely analogous argument, where we reverse the
inequality in the definition of $M^\prime_{i_j}$, shows that $C(p)$ also has a
unique maximal element. Moreover, the same line of reasoning shows that if
$p$ is neither minimal nor maximal in $C(p)$, then there exist $i,j\in
[N]$ such that $M_i(p) \lhd p$ and $M_j(p) \rhd p$.

Let $P^\phi$ be the subposet of $P$ induced by the fixed points of
$\phi$. Any $p\in P^\phi$ is either minimal or maximal
in $C(p)$, because $M_i(p) \lhd p$ either holds for all $i$ or
for none of the $i$. Furthermore, for any given $p\in P$, $\min C(p)$ belongs
to $P^\phi$ if and only if $\max C(p)$ does; this happens
if and only if $\phi(C(p)) = C(p)$.

Assume $p\in P^\phi$ is the maximal element in $C(p)$, and let $q\in P^\phi$
denote the corresponding minimal element. We claim that $p$ covers $q$
in $P^\phi$. Indeed, suppose $r < p$ for some $r\in P^\phi$ with $r =
\min C(r)$. Choose an expression $q = M_{i_t} 
\circ \dots \circ M_{i_1}(p)$ with $t$ as small as possible. Repeated
application of the lifting property shows $M_{i_j} \circ \dots \circ
M_{i_1}(p)\geq r$ for all $j\in [t]$. Hence, $q \geq r$. A similar argument
shows that if $q < r$ and $r = \max C(r)$, then $r \geq p$. The claim is
established.

The above shows that we have a well-defined matching $M^\phi$ on
$P^\phi$ given by
\[
M^\phi(p) = 
\begin{cases}
\min C(p) & \text{if } p = \max C(p),\\
\max C(p) & \text{otherwise.}
\end{cases}
\] 

It remains to show that $M^\phi$ is special, so suppose $p$ covers $q$
in $P^\phi$. First, we assume $p = \max C(p)$. If $\min C(p) = q$,
there is nothing to show. Otherwise, $q = \max C(q)$ as was shown
above. In this case, the above argument shows $\min C(p) > \min C(q)$,
i.e.\ $M^\phi(p) > M^\phi(q)$ as required. The situation when $p = \min
C(p)$ is completely analogous.
\end{proof}

Corollary \ref{co:main} is now straightforward to establish.

\begin{proof}[Proof of Corollary \ref{co:main}]
Let $P$ be a zircon. Choose a non-minimal $p\in P$. The
principal order ideal 
$P_{\leq p} = \{q\in P\mid q\leq p\}$ contains a unique minimal
element $\min P_{\leq p}$. To see
this, choose a special matching $M$ on $P_{\leq p}$ and note that if
this ideal would contain two distinct minimal elements,
then the same would be true for $P_{\leq M(p)}$ because of the lifting
property. Thus, we would obtain
an infinite descending sequence in $P_{\leq p}$ contradicting its
finiteness.

Choose an automorphism $\phi$ of $P$, and let $P^\phi$ denote the
subposet of $P$ induced by the fixed points. If $p\in P^\phi$, then
$\min P_{\leq p}\in P^\phi$, too. Thus, the principal order ideal
$P^\phi_{\leq p}$ coincides with the fixed point set of the 
restriction of $\phi$ to the interval $[\min P_{\leq p}, p]$. By
Theorem \ref{th:main}, $P^\phi_{\leq p}$ has a special matching.
\end{proof}

Earlier, we claimed that our definition of zircons coincides
with the one given by Marietti. We conclude by verifying this
assertion. 

A poset $P$ is {\em ranked} if there is a {\em rank function}
$\rho:P\to \N$ satisfying $\rho(x)=\rho(y)-1$ whenever $x\lhd y$.
 
\begin{definition}[Marietti \cite{marietti}] \label{de:marietti}
A {\em zircon} is a locally finite, ranked poset in which every
non-trivial principal order ideal has a special matching. 
\end{definition}

\begin{proposition}\label{pr:equal}
Definition \ref{de:zircon} and Definition \ref{de:marietti} define the
same class of posets.
\begin{proof}
In this proof, say that a poset satisfying Definition \ref{de:zircon}
is of the {\em first kind}, whereas Definition \ref{de:marietti}
defines posets of the {\em second kind}.

First, suppose $P$ is of the second kind. It was shown in \cite{marietti}
that the principal order ideals in $P$ are intervals, i.e.\ have
unique minimal elements. They are finite since $P$ is locally
finite. Hence, $P$ is of the first kind.

For the converse, we now assume $P$ to be of the first kind. Since every
interval is contained in a principal order ideal, $P$ is locally
finite. To show that $P$ is ranked, it suffices to verify
that for any $p\in P$, the maximal chains in $P_{\leq p}$ all have the
same length. 

Assume, in order to deduce a contradiction, that $P_{\leq p}$ is
a minimal non-ranked principal order ideal. Let $M$ be a special
matching on $P_{\leq p}$. Choose $q\lhd p$ such that the maximal
chains in $P_{\leq q}$ differ in length from those in $P_{\leq
  M(p)}$. Since $M(q)<M(p)$, this implies that there is some $z$ with
$M(q)\lhd z <M(p)$. However, this means that $q \not < M(z)$,
contradicting the fact that $M$ is special. Hence, $P$ is of the
second kind. 
\end{proof}
\end{proposition}

\end{document}